# Bounded nonvanishing functions
# and Bateman functions

Wolfram Koepf, Dieter Schmersau

to appear in: Complex Variables, 1994

**Abstract:**

We consider the family $\widetilde{B}$ of bounded nonvanishing analytic functions $f(z) = a_0 + a_1\,z + a_2\,z^2 + \cdots$ in
the unit disk. The coefficient problem had been extensively investigated ffi(see e. g. [2], [13], [14], [16],
[17], [18], [20]), and it is known that

$$|a_n| \leq \frac{2}{e}$$

for $n = 1, 2, 3$, and $4$. That this inequality may hold for $n \in \mathbb{N}$, is known as the Krzyż conjecture.
It turns out that for $f \in \widetilde{B}$ with $a_0 = e^{-t}$

$$f(z) \prec e^{-t\,\frac{1+z}{1-z}}$$

so that the superordinate functions $e^{-t\,\frac{1+z}{1-z}} = \sum\limits_{k=0}^{\infty} F_k(t)\,z^k$ are of special interest. The corresponding
coefficient functions $F_k(t)$ had been independently considered by Bateman [3] who had ffiintroduced
them with the aid of the integral representation

$$F_k(t) = (-1)^k \frac{2}{\pi} \int\limits_0^{\pi/2} \cos\left(t\,\tan\theta - 2\,k\,\theta\right) d\theta \ .$$

We study the Bateman functions and formulate properties that give insight in the coefficient problem
in $\widetilde{B}$.

## 1  Introduction

We consider functions that are analytic in the unit disk

$$\mathbb{D} := \{z \in \mathbb{C} \mid |z| < 1\} \ .$$

An analytic function $f$ is called subordinate to $g$, if $f = g \circ \omega$ for some analytic function $\omega$ with
$\omega(0) = 0$ and $\omega(\mathbb{D}) \subset \mathbb{D}$; we write $f \prec g$. The subordination principle states that if $g$ is univalent
then $f \prec g$ if and only if $f(0) = g(0)$ and $f(\mathbb{D}) \subset g(\mathbb{D})$, see e. g. [15], § 23.

Let $\widetilde{B}$ denote the family of bounded nonvanishing analytic functions $f(z) = a_0 + a_1\,z + a_2\,z^2 + \cdots$
in $\mathbb{D}$. As $f$ is nonvanishing, we have $\operatorname{Re} \ln f(z) < 0$, and by the subordination principle it turns
out that for $a_0 = e^{-t}$

$$-\ln f(z) \prec t\,\frac{1+z}{1-z} \ ,$$

and so $f(z) = e^{-t} + a_1 z + a_2 z^2 + \cdots \in \widetilde{B}$ if and only if

$$f(z) \prec e^{-t} \frac{1+z}{1-z} \ . \tag{1}$$

Thus the superordinate functions

$$G(t,z) = e^{-t} \frac{1+z}{1-z} =: \sum_{k=0}^{\infty} F_k(t) z^k \tag{2}$$

are of special interest. Graphs of the functions $F_n(t)$ $(n = 1, \ldots, 5)$ are given in Figure 1.

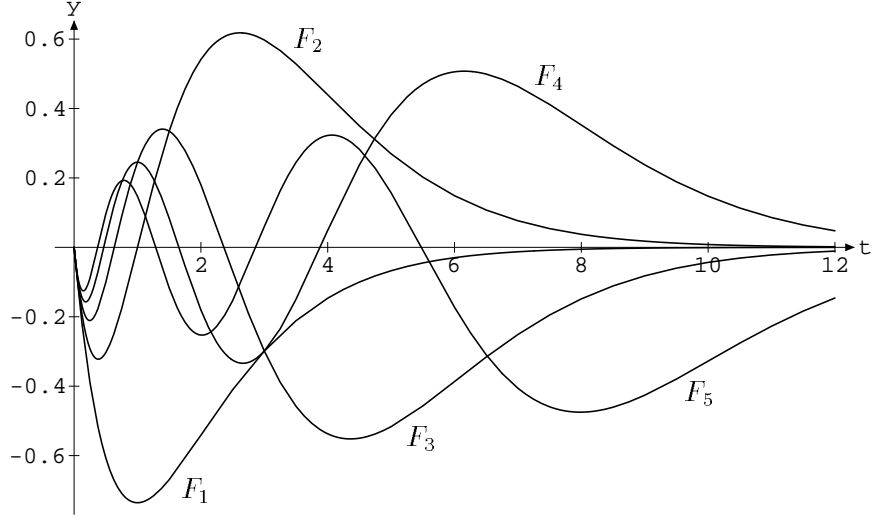

Figure 1: Graphs of the functions $F_n(t)$ $(n = 1, \ldots, 5)$

The following is a list of the first functions $F_n$:

$$
\begin{aligned}
F_0(t) &= e^{-t} \\
F_1(t) &= -2\,t e^{-t} \ , \\
F_2(t) &= 2\,e^{-t}\left(-1 + t\right)t \ , \\
F_3(t) &= \frac{2\,t\,e^{-t}\left(-3 + 6\,t - 2\,t^2\right)}{3} \ , \\
F_4(t) &= \frac{2\,t\,e^{-t}\left(-3 + 9\,t - 6\,t^2 + t^3\right)}{3} \ , \\
F_5(t) &= \frac{2\,t\,e^{-t}\left(-15 + 60\,t - 60\,t^2 + 20\,t^3 - 2\,t^4\right)}{15} \ .
\end{aligned}
$$

We consider the coefficient problem $(n \in \mathbb{N})$ to find

$$A_n := \max_{f \in \widetilde{B}} |a_n(f)| \ .$$

That the maximum exists for all $n \in \mathbb{N}$ follows from the fact that the union of $\widetilde{B}$ with the constant functions $c \in \overline{\mathbb{D}}$ forms a compact family of analytic functions. For the coefficient problem it is no loss of generality to assume that $a_0 > 0$ so that we can assume that (1) holds for some $t > 0$. For small $n$ it is then easy to solve the coefficient problem using subordination techniques.

As $f \prec g$ implies that $|a_1(f)| \leq |a_1(g)|$ (see e. g. [15], Theorem 212), we have

$$A_1 = \max_{t \geq 0} \left| a_1 \left( e^{-t \frac{1+z}{1-z}} \right) \right| = \max_{t \geq 0} |F_1(t)| = \max_{t \geq 0} 2\, t\, e^{-t} = \frac{2}{e}$$

with equality iff $t = 1$, and $f(z) = \eta\, e^{-\frac{1+\xi z}{1-\xi z}}$ ($|\xi| = |\eta| = 1$).

By the composition with a Möbius transform, this leads to the inequality (see [16])

$$(1 - |z|^2)\, |f'(z)| \leq \frac{2}{e} \qquad (z \in \mathbb{D}) \tag{3}$$

from which we may deduce by a standard technique (see e. g. [7], p. 72, Exercise 17) that

$$
\begin{aligned}
n\, |a_n(f)| &= |a_{n-1}(f')| = \left| \frac{1}{2\pi i} \int_{\partial \mathbb{D}_r} \frac{f'(\zeta)}{\zeta^n} d\zeta \right| \leq \frac{1}{r^{n-1}} \left( \frac{1}{2\pi} \int_0^{2\pi} |f'(re^{i\theta})|\, d\theta \right) \\
&\leq \frac{1}{r^{n-1}(1 - r^2)} \frac{2}{e} \leq \frac{2}{e} \frac{n+1}{2} \left( 1 + \frac{2}{n-1} \right)^{\frac{n-1}{2}} < n
\end{aligned}
$$

where we used (3) and chose $r^2 = \frac{n-1}{n+1}$. Unfortunately this estimate is just too weak to be offff value: The bound 1 for $A_n$ is very elementary, and holds even for all functions bounded by 1ff. Each global bound less than 1 would be new, however.

It is similarly easy to solve the coefficient problem for $n = 2$ using subordination techniques, ([15], Theorem 212, see e. g. [10]).

Using several methods it was shown that

$$A_n = \frac{2}{e} \tag{4}$$

for $n = 1, 2, 3,$ and 4. Obviously $G(1, z^n)$ has $n^{th}$ coefficient equal to $\frac{2}{e}$, which makes these results sharp. That (4) may hold for $n \in \mathbb{N}$, is known as the Krzyż conjecture.

If the subordinate function has very regular coefficient behffavior, then global coefficient results are available: If

$$\sum_{k=0}^{\infty} a_k\, z^k \prec \sum_{k=0}^{\infty} b_k\, z^k\ , \tag{5}$$

and if the coefficient sequence $b_n$ is nonnegative, nonincreasing, and convex, then $|a_n| \leq b_0$ for all $n \in \mathbb{N}_0$, and if the coefficient sequence $b_n$ is nonnegative, nondecreasing, and convex, then $|a_n| \leq b_n$ all $n \in \mathbb{N}_0$ (see e.g. [15], Theorem 216). On the other hand, the coefficient sequencesff of our subordinate functions $G(t, z)$ are highly irregular for all $t > 0$.

Another important result, however, can be obtained by subordination techniques, as well. It is well known that if (5) holds, then

$$\sum_{k=0}^{\infty} |a_k|^2 \leq \sum_{k=0}^{\infty} |b_k|^2\ .$$

Especially: If an analytic function $f(z) = \sum\limits_{k=0}^{\infty} a_k \, z^k$ of the unit disk is bounded by 1, then $f \prec z$, and the relation

$$\sum_{k=0}^{\infty} |a_k|^2 \leq 1$$

(following also directly from Parseval's identity) is obtffained (for sharper versions see also [7], Theorem 6.1 and Theorem 6.2). Equality holds if and only if (see [19], Theorem 3) $f$ is an inner function, i.e. if the radial limit $f(e^{i\theta}) := \lim\limits_{r \to 1} f(r e^{i\theta}) = 1$ for almost all $e^{i\theta}$ on the unit circle $\partial \mathbb{D}$. Nonvanishing inner functions with positive $f(0)$ have the representation (see e.g. [9], second theorem ffon p. 66)

$$f(z) = \exp\left( - \int\limits_{\partial \mathbb{D}} \frac{1 + e^{i\theta} z}{1 - e^{i\theta} z} \, d\mu(\theta) \right),$$

where $\mu$ is a singular positive measure on the unit circle $\partial \mathbb{D}$. If we choose a point measure $\mu$ concentrating its full measure $t$ at the point $\theta = 1$, we get the function $G(t,z) = e^{-t \frac{1+z}{1-z}}$ of Equation (2) so that we are lead to the identity

$$\sum_{k=0}^{\infty} F_k^2(t) = 1 \; .$$

For each individual coefficient of $G(t,z)$ we thus have the (weak) inequality

$$|F_n(t)| \leq \sqrt{1 - F_0^2(t)} = \sqrt{1 - e^{-2t}} \qquad (n \in \mathbb{N}) \; . \tag{6}$$

It is the purpose of this paper to develop further properties, especially inequalities, for the functions $F_n$ $(n \in \mathbb{N})$, giving more insight in the coefficient problem for $\widetilde{B}$.

## 2   A collection of properties of the Bateman functions

In [3] (see also [1], § 13.6) Bateman introduced the functions $(x \geq 0)$

$$k_n(x) := \frac{2}{\pi} \int\limits_0^{\pi/2} \cos\left( x \, \tan\theta - n \, \theta \right) d\theta \; ,$$

and he verified that ([3], formula (2.7))

$$k_{2m}(x) = (-1)^m \, e^{-x} \Big( L_m(2x) - L_{m-1}(2x) \Big) \tag{7}$$

where $L_m(t)$ denotes the $m^{th}$ Laguerre polynomial. On the other hand if one defines the functions $F_n$ $(n \in \mathbb{N}_0)$ with the aid of the generating function

$$e^{-t \frac{1+z}{1-z}} =: \sum_{k=0}^{\infty} F_k(t) \, z^k \; ,$$

one gets immediately (see [10], formula (14), and p. 178)

$$F_n(t) = e^{-t} \left( L_n(2t) - L_{n-1}(2t) \right) \tag{8}$$

and a comparison of (7) and (8) yields the relation

$$F_n(t) = (-1)^n \, k_{2n}(t)$$

so that we get the Bateman representation

$$F_n(t) = (-1)^n \frac{2}{\pi} \int\limits_0^{\pi/2} \cos\left( t \tan\theta - 2\, n\, \theta \right) d\theta \tag{9}$$

for our functions $F_n$. By Bateman's work we are prepared to state many further properties: For $n \in \mathbb{N}$ the function $F_n$ satisfies the differential equation (see [3], formula (5.1))

$$t\, F_n''(t) = (t - 2n)\, F_n(t) \tag{10}$$

with the initial values

$$F_n(0) = 0 \qquad \text{and} \qquad F_n'(0) = -2 \,, \tag{11}$$

and the Rodriguez type formula (see [3], formula (31))

$$F_n(t) = \frac{t\, e^t}{n!} \frac{d^n}{dt^n} \left( e^{-2t}\, t^{n-1} \right) \,.$$

The differential equation can also be obtained completely algorithmically (see [11]–[12]).

Further we get the following connection wiffth the generalized Laguerre polynomials (see [23], ffp. 216, formula (1.15))

$$F_n(t) = e^{-t}\, L_n^{(-1)}(2t) \,, \tag{12}$$

and (see [22], formula (5.2.1))

$$F_n(t) = -e^{-t} \frac{2t}{n}\, L_{n-1}^{(1)}(2t) \,, \tag{13}$$

from which one may deduce the hypergeometric representation

$$F_n(t) = -2t\, e^{-t} {}_1F_1 \left( \begin{array}{c} 1 - n \\ 2 \end{array} \middle|\; 2t \right) \,,$$

and the explicit representation

$$F_n(t) = \frac{e^{-t}}{n} \sum_{k=1}^n \frac{(-1)^k}{(k-1)!} \binom{n}{k} (2t)^k \,.$$

Bateman obtained further relations: a difference equation ([3], formula (4.1))

$$(n-1)(F_n(t) - F_{n-1}(t)) + (n+1)(F_n(t) - F_{n+1}(t)) = 2\, t\, F_n(t) \tag{14}$$

that is also an easy consequence of the defining equation using the generating functioffn, he obtained a difference differential equation ([3], formula (4.2))

$$(n+1) F_{n+1}(t) - (n-1) F_{n-1}(t) = 2\,t\,F_n'(t)\,,$$

and a system of differential equations ([3], formula (4.3))

$$F_n'(t) - F_{n+1}'(t) = F_n(t) + F_{n+1}(t)\,, \tag{15}$$

from which he is lead to the inequalities for $F_n$ ([3], formula (4.4))

$$|F_n(t)| \leq \frac{2\,n}{t} \quad (n > 2)\,, \tag{16}$$

and for $F_n'$ ([3], formula (4.5))

$$|F_n'(t)| \leq \frac{n}{t} \quad (n > 2)\,. \tag{17}$$

For large $t$ the first inequality is a refinement of the trivial estimate

$$|F_n(t)| \leq 1 \tag{18}$$

that follows from (6) or from the Bateman representation (9).

Finally Bateman obtained the following statements about Integrals of products $F_n\,F_m$ ($n, m \in \mathbb{N}$) (see [3], formula (2.91)

$$\int\limits_0^\infty F_n^2(t)\,dt = 1 \qquad \text{and} \qquad \int\limits_0^\infty F_n(t)\,F_m(t)\,dt = \begin{cases} 0 & \text{if } |n-m| \neq 1 \\ \frac{1}{2} & \text{otherwise} \end{cases}\,. \tag{19}$$

We state further properties: The functions $F_n$ ($n \in \mathbb{N}$) have a zero at the origin and $n-1$ further positive real zeros (see e. g. [23], *Nullstellensatz*, p. 123) (indeed, by (13), $F_n(t)$ has the same zeros as $L_{n-1}^{(1)}(2t)$).

From the differential equation (10) we moreover see that at $t = 2n$ there is a point of inflection, and as $F_n(t) \to 0$ for $t \to \infty$, and all other points of inflection lie at the zeros of $F_n$ one easily deduces that $t = 2n$ must be the largest point of inflection of $F_n$ implying that all the zeros of $F_n$ lie in the interval $[0, 2n]$. The successive relative maxima of $|F_n|$ lying between the zeros of $F_n$ form an increasing sequence (see [22], Theorem 7.6.2, $\alpha = -1$), so that the largest value attained by $|F_n(t)|$ is attained at the last zero of $F_n'$ which is seen to lie between the last zero $T_n$ of $F_n$ and the point $t = 2n$. For small $n$ the mentioned qualitative properties of $F_n$ can be recognized in Figure 1.

By a result of Hahn ([8], formula (17)) the last zero $T_n$ (being the last zero of $L_{n-1}^{(1)}(2t)$) satisfies the relation

$$4n - 2 - C_1 \sqrt[3]{4n-2} < 2\,T_n < 4n - 2 - C_2 \sqrt[3]{4n-2} \tag{20}$$

with two positive constants $C_1, C_2 \in \mathbb{R}^+$ that are independent of $n$, in particular

$$\lim_{n \to \infty} \frac{T_n}{n} = 2\,. \tag{21}$$

The right hand side of (20) leads to the sharpened inequalitffy

$$T_n < 2n - 1\,,$$

and Puiseux series expansion of (20) yields the refinement of (21)

$$\frac{T_n}{n} = 2 - O\left(\left(\frac{1}{n}\right)^{2/3}\right)\,.$$

# 3   Representation by residues

To the system of differential equations given by (15) together with the initial conditions $F_n(0) = 0$ ($n \in \mathbb{N}$) the technique of Laplace transformation

$$\mathcal{L}(f)(z) := \int\limits_0^\infty e^{-zt} f(t) \, dt$$

can be applied to deduce a representation by residues for $F_n$. It is well-known that $\mathcal{L}(f') = z\,\mathcal{L}(f) - f(0)$ (see e.g. [6], Satz 9.1) so that we obtain ($n \in \mathbb{N}$)

$$(z+1)\,\mathcal{L}(F_{n+1}) = (z-1)\,\mathcal{L}(F_n)$$

or

$$\mathcal{L}(F_{n+1}) = \frac{z-1}{z+1}\,\mathcal{L}(F_n) \, .$$

Induction shows then that for $n \in \mathbb{N}$ and $k \in \mathbb{N}_0$

$$\mathcal{L}(F_{n+k}) = \left(\frac{z-1}{z+1}\right)^k \mathcal{L}(F_n) \, . \tag{22}$$

To obtain the initial function $\mathcal{L}(F_1)$, we use $F_0(t) = e^{-t}$ to get first

$$\mathcal{L}(F_0)(z) = \int\limits_0^\infty e^{(z-1)t} \, dt = \frac{1}{1+z} \, .$$

Further from (15) with $n = 0$ we are lead to

$$(z+1)\,\mathcal{L}(F_1) = (z-1)\,\mathcal{L}(F_0) - 1$$

or

$$\mathcal{L}(F_1)(z) = -\frac{2}{(1+z)^2} \, .$$

Thus by an application of (22) with $n = 1$ we have finally

$$\mathcal{L}(F_k)(z) = -\frac{2}{(1+z)^2} \left(\frac{z-1}{z+1}\right)^{k-1} \, .$$

If we use now the inverse Laplace transform (see e.g. [6], p. 170, fformula (15)), we get

$$F_k(t) = \lim_{R \to \infty} \frac{1}{2\pi i} \int\limits_{\gamma_R} e^{tz} \mathcal{L}(F_k)(z) \, dz \, ,$$

where $\gamma_R : [-R, R] \to \mathbb{C}$ is given by $\gamma_R(\tau) = i\tau$, and therefore we have the integral representation

$$F_k(t) = \frac{1}{\pi} \int\limits_{-\infty}^\infty e^{it\tau} \frac{(\tau+i)^{k-1}}{(\tau-i)^{k+1}} \, d\tau \, .$$

By a standard procedure this can be identified as the residue (see e.g. [4ff], p. 217, formula (12))

$$\int\limits_{-\infty}^{\infty} e^{it\tau} \frac{(\tau + i)^{k-1}}{(\tau - i)^{k+1}}\, d\tau = 2\pi i \operatorname{Res}\left(e^{itz}\frac{(z + i)^{k-1}}{(z - i)^{k+1}}\right)\ ,$$

and therefore we have the representation ($k \in \mathbb{N}$)

$$F_k(t) = 2i \operatorname{Res}\left(e^{itz}\frac{(z + i)^{k-1}}{(z - i)^{k+1}}\right)\ .$$

These results are collected in

**Theorem 1** The Bateman functions $F_k$ ($k \in \mathbb{N}$) satisfy the integral representation

$$F_k(t) = \frac{1}{\pi} \int\limits_{-\infty}^{\infty} e^{it\tau} \frac{(\tau + i)^{k-1}}{(\tau - i)^{k+1}}\, d\tau$$

and therefore the residual representation

$$F_k(t) = 2i \operatorname{Res}\left(e^{itz}\frac{(z + i)^{k-1}}{(z - i)^{k+1}}\right)\ .$$

## 4  Results deduced from the differential equation

In this section we deduce another statements about an integral involving the Bateman functions and get an estimate for $|F_n|$ using its differential equation (10). Multiplying (10) by $2\, F_n'(t)/t$, we have

$$2\, F_n'(t)\, F_n''(t) = 2\, F_n(t)\, F_n'(t) - \frac{2n}{t}\, 2\, F_n(t)\, F_n'(t)\ .$$

We integrate from $0$ to $t$, and get for $n \in \mathbb{N}$ using the initial values (11)

$$\left(F_n'\right)^2(t) - 4 = F_n^2(t) - \int\limits_0^t \frac{2n}{\tau}\, 2\, F_n(\tau)\, F_n'(\tau)\, d\tau\ . \tag{23}$$

For the last integral we get integrating by parts

$$
\begin{aligned}
\int\limits_0^t \frac{2n}{\tau}\, 2\, F_n(\tau)\, F_n'(\tau)\, d\tau 
&= \left.\frac{2n}{\tau} F_n^2(\tau)\right|_0^t + \int\limits_0^t \frac{2n}{\tau^2} F_n^2(\tau)\, d\tau \\[2mm]
&= \frac{2n}{t}\, F_n^2(t) - 2n\, F_n(0)\, F_n'(0) + \int\limits_0^t 2n\, \left(\frac{F_n(\tau)}{\tau}\right)^2 d\tau \\[2mm]
&= \frac{2n}{t}\, F_n^2(t) + 2n \int\limits_0^t \left(\frac{F_n(\tau)}{\tau}\right)^2 d\tau\ .
\end{aligned}
$$

So we have the identity

$$(F_n')^2 (t) - 4 = F_n^2(t) - \frac{2n}{t} F_n^2(t) - 2n \int_0^t \left( \frac{F_n(\tau)}{\tau} \right)^2 d\tau \ . \tag{24}$$

From this identity by letting $t \to \infty$ we are lead to the statement

$$\int_0^\infty \left( \frac{F_n(\tau)}{\tau} \right)^2 d\tau = \frac{2}{n}$$

as $\lim_{t \to \infty} F_n(t) = \lim_{t \to \infty} F_n'(t) = 0$. Therefore in particular $(t \geq 0)$

$$\int_0^t \left( \frac{F_n(\tau)}{\tau} \right)^2 d\tau < \frac{2}{n} \ . \tag{25}$$

At this point we like to mention that from (24) it is now very easy to deduce the inequality $t < 2n$ for a local extremum of $F_n$, again (compare § 2), as an application of (25) yields

$$F_n^2(t) - \frac{2n}{t} F_n^2(t) = (F_n')^2 (t) - 4 + 2n \int_0^t \left( \frac{F_n(\tau)}{\tau} \right)^2 d\tau < (F_n')^2 (t) \ ,$$

and therefore for any point with $F_n'(t) = 0$ we get $t < 2n$.

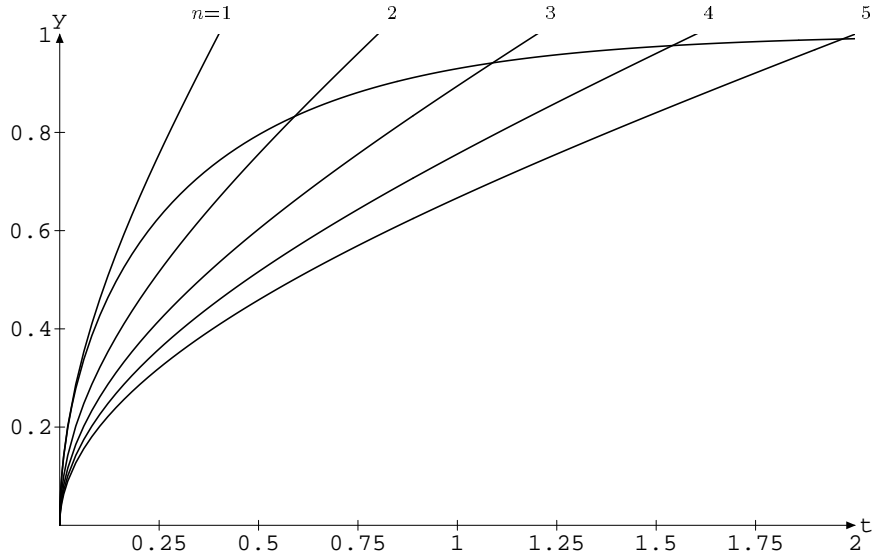

Figure 2: The estimates (6) and (26) for $n = 1, \ldots, 5$

To deduce an estimate for $|F_n|$ we regroup (24) and get

$$\frac{2n - t}{t} F_n^2(t) = 4 - (F_n')^2 (t) - 2n \int_0^t \left( \frac{F_n(\tau)}{\tau} \right)^2 d\tau < 4 \ ,$$

and for $t < 2n$ (which is the critical region) finally

$$|F_n^2(t)| < \frac{4t}{2n - t} \ . \tag{26}$$

We note that, however, this improves (6) for small $t$ only, see Figure 2.

In the next section we will give a further improvemeffnt of (26).

## 5  Estimates by the Szegö method

We consider the generalized Laguerre polynomials $L_n^{(\alpha)}(x)$ $(\alpha \in \mathbb{R})$ given by their Rodriguez formula

$$x^\alpha \, e^{-x} \, L_n^{(\alpha)}(x) = \frac{1}{n!} \frac{d^n}{dx^n} \left( e^{-x} \, x^{n+\alpha} \right) \tag{27}$$

(see [23], p. 213, formula (1.3)). Szegffff7ö ([21], see [22], p. 159, formula (7.21.3)) considered the case $\alpha = 0$, and was lead to the inequality

$$e^{-\frac{x}{2}} \left| L_n^{(0)}(x) \right| < 1 \qquad (x > 0) \ . \tag{28}$$

Using a similar method we get the following development. For $n \in \mathbb{N}_0$ and $\alpha \in \mathbb{Z}$ the function $f_{n\alpha}(x) := e^{-x} \, x^{n+\alpha}$ is analytic in $\{z \in \mathbb{C} \mid \operatorname{Re} z > 0\}$. By (27) we have

$$x^\alpha \, e^{-x} \, L_n^{(\alpha)}(x) = \frac{1}{n!} f_{n\alpha}^{(n)}(x) \ .$$

If $x \in \mathbb{R}^+$, then for $z = x + re^{i\theta}$ $(r \in (0, x))$ we have by Taylor's formula

$$f_{n\alpha}(z) = \sum_{k=0}^{\infty} \frac{1}{k!} f_{n\alpha}^{(k)}(x) \, (z - x)^k \ ,$$

and Cauchy's integral formula gives the estimate

$$\frac{1}{k!} \left| f_{n\alpha}^{(k)}(x) \right| < \frac{1}{r^k} \max_{0 \le \theta \le 2\pi} |f_{n\alpha}(x + re^{i\theta})| \ .$$

Especially for $k = n$ this yields

$$x^\alpha \, e^{-x} \left| L_n^{(\alpha)}(x) \right| r^n < \max_{0 \le \theta \le 2\pi} |f_{n\alpha}(x + re^{i\theta})| \tag{29}$$

If furthermore $n + \alpha > 0$, then $f_{n\alpha}$ is analytic in all of $\mathbb{C}$, and (29) holds for all $r \in \mathbb{R}^+$. This case will be studied now.

To give an estimate of $\max_{0 \le \theta \le 2\pi} |f_{n\alpha}(x + re^{i\theta})|$, we expand

$$f_{n\alpha}(x + re^{i\theta}) = e^{-(x + re^{i\theta})} \, (x + re^{i\theta})^{n+\alpha}$$

to get

$$\max_{0 \le \theta \le 2\pi} |f_{n\alpha}(x + re^{i\theta})| = e^{-x} \max_{0 \le \theta \le 2\pi} e^{-r\cos\theta} \, (x^2 + r^2 + 2rx\cos\theta)^{\frac{n+\alpha}{2}} \ .$$

Together with (29) we have therefore

$$x^\alpha \left| L_n^{(\alpha)}(x) \right| r^n < \max_{0 \le \theta \le 2\pi} e^{-r\cos\theta} \left( x^2 + r^2 + 2rx\cos\theta \right)^{\frac{n+\alpha}{2}} \, .$$

We set now $\lambda := \cos\theta \in [-1, 1]$, $p(\lambda) := x^2 + r^2 + 2rx\lambda$, and $q(\lambda) := e^{-r\lambda} p(\lambda)^{\frac{n+\alpha}{2}}$, and have therefore

$$x^\alpha \left| L_n^{(\alpha)}(x) \right| r^n < \max_{-1 \le \lambda \le 1} q(\lambda) \, . \tag{30}$$

As

$$q'(\lambda) = r \, e^{-r\lambda} \, p(\lambda)^{\frac{n+\alpha-2}{2}} \left( x(n+\alpha) - p(\lambda) \right)$$

we get for a possible critical point $\lambda_0$ of $q$ the relation

$$\lambda_0 = \frac{x(n+\alpha) - (x^2 + r^2)}{2rx} \, .$$

At the point $\lambda = \lambda_0$ we have furthermore

$$p(\lambda_0) = x^2 + r^2 + 2rx\lambda_0 = x(n+\alpha) > 0 \, ,$$

hence

$$q''(\lambda_0) = -2r^2 \, x \, e^{-r\lambda_0} \, p(\lambda_0)^{\frac{n+\alpha-2}{2}} < 0 \, ,$$

and $\lambda_0$ maximizes $q$. Therefore, from (30) we get

$$x^\alpha \left| L_n^{(\alpha)}(x) \right| r^n < q(\lambda_0) \tag{31}$$

if $-1 \le \lambda_0 \le 1$.

We consider now the case $x \in (0, 4(n+\alpha)]$ (which with respect to the representation (13) corresponds to the critical region $(0, 2n)$ for $t = x/2$), and choose $r := \sqrt{x(n+\alpha)}$. In this case we have $\lambda_0 = \frac{1}{2}\sqrt{\frac{x}{n+\alpha}} \in [-1, 0)$. Hence, (31) implies

$$x^\alpha \left| L_n^{(\alpha)}(x) \right| (x(n+\alpha))^{\frac{n}{2}} < e^{\frac{x}{2}} \left( x(n+\alpha) \right)^{\frac{n+\alpha}{2}} \, ,$$

and therefore finally

**Theorem 2** For the generalized Laguerre polynomials $L_n^{(\alpha)}(x)$ ($\alpha \in \mathbb{Z}$) the estimate

$$e^{-\frac{x}{2}} \left| L_n^{(\alpha)}(x) \right| < \left( \frac{n+\alpha}{x} \right)^{\frac{\alpha}{2}} \tag{32}$$

holds for $x \in (0, 4(n+\alpha)]$ if $n + \alpha > 0$. $\qquad\square$

If we define the functions ($\alpha \in \mathbb{R}$)

$$F_n^{(\alpha)}(t) := e^{-t} L_n^{(\alpha)}(2t) \tag{33}$$

then (32) reads ($x = 2t$)

$$\left| F_n^{(\alpha)}(t) \right| < \left( \frac{n+\alpha}{2t} \right)^{\frac{\alpha}{2}} \qquad \left( t \in (0, 2(n+\alpha)] \right) \, .$$

For $\alpha = 0$, we have Szegö's result (28) in this interval, and for $\alpha = 1$ we get in view of representation (13)

$$|F_n(t)| < \sqrt{\frac{2t}{n}} \qquad (t \in (0, 2n)) .$$

(34)

This inequality improves (26) as a Puiseux expansion yields

$$\sqrt{\frac{4t}{2n-t}} = \sqrt{\frac{2t}{n}} + \frac{1}{2\sqrt{2}} \left( \frac{t}{n} \right)^{\frac{3}{2}} + P \left( \frac{t}{n} \right)$$

with some positive function $P$.

Note that the special choice $\alpha = -1$ (and *not* the value $\alpha = 0$) generates the Bateman functions $F_n(t) = F_n^{(-1)}(t)$.

# 6 Asymptotic estimates

We consider the functions ($\alpha \in \mathbb{R}$)

$$F_n^{(\alpha)}(t) = e^{-t} L_n^{(\alpha)}(2t) = e^{-t} \left( L_n^{(\alpha+1)}(2t) - L_{n-1}^{(\alpha+1)}(2t) \right) = F_n^{(\alpha+1)}(t) - F_{n-1}^{(\alpha+1)}(t)$$

(35)

of (33) (see [23], p. 216, formula (1.15)) now in more detail. Taking derivative yields

$$
\begin{aligned}
\left( F_n^{(\alpha)} \right)'(t) &= -F_n^{(\alpha)}(t) + 2\, e^{-t} \left( L_n^{(\alpha)} \right)'(2t) \\
&\overset{(35)}{=} -e^{-t} \left( L_n^{(\alpha+1)}(2t) - L_{n-1}^{(\alpha+1)}(2t) \right) + 2\, e^{-t} \left( L_n^{(\alpha)} \right)'(2t) \\
&= -e^{-t} \left( L_n^{(\alpha+1)}(2t) + L_{n-1}^{(\alpha+1)}(2t) \right) \\
&= -\left( F_n^{(\alpha+1)}(t) + F_{n-1}^{(\alpha+1)}(t) \right) .
\end{aligned}
$$

(36)

where the relation about $\left( L_n^{(\alpha)} \right)'$ corresponds to ([23], p. 215, formula (1.12)).

Moreover the program [12] generates the differential equation

$$(1 + \alpha + 2\, n - t)\, F_n^{(\alpha)}(t) + (1 + \alpha) \left( F_n^{(\alpha)} \right)'(t) + t \left( F_n^{(\alpha)} \right)''(t) = 0$$

for the functions $F_n^{(\alpha)}$ with respect to the variable $t$, and the recurrence equation

$$(-1 + \alpha + n)\, F_{n-2}^{(\alpha)} + (1 - \alpha - 2\, n + 2\, t)\, F_{n-1}^{(\alpha)} + n\, F_n^{(\alpha)} = 0$$

with respect to the variable $n$ (check!).

Assume now, $t_1 < t_2$. Then we get by an integration

$$
\begin{aligned}
\left( F_n^{(\alpha)} \right)^2 (t_2) - \left( F_n^{(\alpha)} \right)^2 (t_1) &= 2 \int_{t_1}^{t_2} F_n^{(\alpha)}(t) \left( F_n^{(\alpha)} \right)'(t)\, dt \\
&\overset{(35),(36)}{=} -2 \int_{t_1}^{t_2} \left( \left( F_n^{(\alpha+1)}(t) \right)^2 - \left( F_{n-1}^{(\alpha+1)}(t) \right)^2 \right) dt .
\end{aligned}
$$

If we choose $\alpha = 0$, we get in particular

$$\left(F_n^{(0)}\right)^2(t_2) - \left(F_n^{(0)}\right)^2(t_1) = -2 \int_{t_1}^{t_2} \left(\left(F_n^{(1)}(t)\right)^2 - \left(F_{n-1}^{(1)}(t)\right)^2\right) dt \ . \tag{37}$$

Together with the relation ($n \in \mathbb{N}$)

$$F_n(t) = F_n^{(-1)}(t) = -\frac{2t}{n} F_{n-1}^{(1)}(t)$$

(see (12) and (13), or [22], p. 98, formula (5.2.1)) it followffs from (37) that

$$\left(F_n^{(0)}\right)^2(t_2) - \left(F_n^{(0)}\right)^2(t_1) = \int_{t_1}^{t_2} F_{n+1}^2(t)\frac{(n+1)^2}{2\,t^2}\,dt - \int_{t_1}^{t_2} F_n^2(t)\frac{n^2}{2\,t^2}\,dt \ .$$

We now let $t_2 \to \infty$. Then $F_n^{(0)}(t_2) \to 0$, and as

$$\int_0^\infty F_n^2(t)\,dt = 1 \tag{38}$$

(see (19), compare [3], formulae (2.7) and (2.91)), we get further

$$\begin{aligned}
\left(F_n^{(0)}\right)^2(t_1) &= \int_{t_1}^\infty F_{n+1}^2(t)\frac{(n+1)^2}{2\,t^2}\,dt - \int_{t_1}^\infty F_n^2(t)\frac{n^2}{2\,t^2}\,dt \\
&\leq \int_{t_1}^\infty F_{n+1}^2(t)\frac{(n+1)^2}{2\,t^2}\,dt \leq \frac{(n+1)^2}{2\,t_1^2}\left(1 - \int_0^{t_1} F_{n+1}^2(t)\,dt\right)
\end{aligned}$$

where the last inequality is deduced with (38). Soff we have finally the inequality ($n \in \mathbb{N}$, $t > 0$)

$$\left|F_n^{(0)}(t)\right| \leq \frac{n+1}{\sqrt{2}\,t} \ .$$

This sharpens the result of Szegö (28) for large $t$. From (35) and (36) it follows that

$$F_n^{(\alpha)}(t) + \left(F_n^{(\alpha)}\right)'(t) = -2\,F_{n-1}^{(\alpha+1)}(t) \ .$$

We deduce now for a critical point $t_k$ of $F_n$ with $F_n'(t_k) = 0$ the relation

$$|F_n(t_k)| = 2\left|F_{n-1}^{(0)}(t_k)\right| \leq \sqrt{2}\frac{n}{t_k} \ . \tag{39}$$

Especially is this $\leq 2/e$ for

$$\frac{t_k}{n} \geq \frac{e}{\sqrt{2}} \approx 1.92211551407955841... \ . \tag{40}$$

It is now remarkable that by the result of Hahn (21) for $n \to \infty$ the most important critical point $T$ of $F_n$ which produces the maximal value of $F_n$ has the property $T/n \to 2$ as $T_n < T < 2$. This gives finally the following

**Theorem 3** The Krzyż conjecture is asymptotically true for the superordinate functions $e^{-t\frac{1+z}{1-z}}$, i.e. we have $\left| a_n\left(e^{-t\frac{1+z}{1-z}}\right)\right| \leq \frac{2}{e}$ for $n \geq N$. $\qquad\qquad\square$

We will now strengthen this result.

Therefore let an arbitrary positive zero $t_n$ of $F_n$ be given. Then $t_n$ is also a zero of $h_n := F_n^2$, and as

$$h_n''(t) = 2\left(F_n'(t)\right)^2 + 2\,F_n(t)\,F_n''(t)\,,$$

by the differential equation for $F_n$ we get

$$h_n''(t) = 2\left(F_n'(t)\right)^2 + 2\,\frac{t-2n}{t}\,F_n^2(t)\,. \tag{41}$$

From this we may deduce that $h_n''(t) > 0$ for $t \geq 2n$. Now, however, we consider the interval between $t_n$ and the smallest relative extremum $t_n^* > t_n$ of $F_n$, i. e. the smallest zero $t_n^*$ of $F_n'$ after $t_n$. Then obviously $h_n$ is strictly increasing in $[t_n, t_n^*]$, further $h_n'(t_n) = h_n'(t_n^*) = 0$, and therefore $h_n'\big|_{[t_n, t_n^*]}$ assumes an absolute maximum at some interior point $t_n^{**} \in (t_n, t_n^*)$, where $h_n''(t_n^{**}) = 0$.

From (41) we deduce

$$\left(F_n'(t_n^{**})\right)^2 = \frac{2n - t_n^{**}}{t_n^{**}}\,F_n^2(t_n^{**})\,,$$

and therefore by (18)

$$\left|F_n'(t_n^{**})\right| = \sqrt{\frac{2n - t_n^{**}}{t_n^{**}}}\,\left|F_n(t_n^{**})\right| < \sqrt{\frac{2n - t_n^{**}}{t_n^{**}}}\,.$$

As $\sqrt{\frac{2n-t}{t}}$ is strictly decreasing for $t \in (0, 2n)$, it follows furthermore that

$$\left|F_n'(t_n^{**})\right| < \sqrt{\frac{2n - t_n}{t_n}}\,,$$

and finally ($h_n'$ is positive)

$$h_n'(t_n^{**}) = 2\left|F_n(t_n^{**})\right|\left|F_n'(t_n^{**})\right| < 2\sqrt{\frac{2n - t_n}{t_n}}$$

using (18) again. As $t_n^{**}$ is the global maximum of $h_n'$ in $[t_n, t_n^*]$, we therefore are lead to the inequalities

$$0 < h_n'(t) < 2\sqrt{\frac{2n - t_n}{t_n}} \tag{42}$$

for all $t \in (t_n, t_n^*)$.

We are interested in $h_n(t_n^*)$, the value of $h_n$ at its maximum $t_n^*$. Therefore let $p > 0$ be given such that $h_n(t_n^*) > \frac{1}{p}$. As $h_n(t_n) = 0$, and $h_n$ is stricly increasing, there is some $\widetilde{t_n} \in (t_n, t_n^*)$ with $h_n(\widetilde{t_n}) = \frac{1}{p}$. The mean value theorem then shows the existence of $\xi_n \in (\widetilde{t_n}, t_n^*)$ with

$$\frac{h_n(t_n^*) - h_n(\widetilde{t_n})}{t_n^* - \widetilde{t_n}} = h_n'(\xi_n)\,,$$

and therefore by (42)

$$\frac{h_n(t_n^*) - h_n(\widetilde{t_n})}{t_n^* - \widetilde{t_n}} < 2 \sqrt{\frac{2n - t_n}{t_n}}$$

or

$$h_n(t_n^*) < h_n(\widetilde{t_n}) + 2 \sqrt{\frac{2n - t_n}{t_n}} (t_n^* - \widetilde{t_n}) = \frac{1}{p} + 2 \sqrt{\frac{2n - t_n}{t_n}} (t_n^* - \widetilde{t_n}) .$$

By (19) we have

$$\int_0^\infty h_n(\tau) \, d\tau = 1 ,$$

and thus by the integral mean value theorem ($h_n$ is nonnegative)

$$1 > \int_{\widetilde{t_n}}^{t_n^*} h_n(\tau) \, d\tau = h_n(\eta_n) (t_n^* - \widetilde{t_n})$$

for some $\eta_n \in (\widetilde{t_n}, t_n^*)$. As $h_n$ is stricly increasing, we therefore get $h_n(\eta_n) > h_n(\widetilde{t_n}) = \frac{1}{p}$ implying

$$1 > \frac{1}{p} (t_n^* - \widetilde{t_n}) \qquad \text{or} \qquad t_n^* - \widetilde{t_n} < p .$$

Finally we have

$$h_n(t_n^*) < \frac{1}{p} + 2p \sqrt{\frac{2n - t_n}{t_n}} .$$

We were lead to this inequality under the assumption that $h_n(t_n^*) > \frac{1}{p}$. If $h_n(t_n^*) \leq \frac{1}{p}$, however, then the same conclusion follows trivially, so that the above calculations can be summarized by the following

**Lemma 1** Let $h_n(t) = F_n^2(t)$, let $t_n$ be a positive zero of $F_n$, let $t_n^*$ the lowest zero of $F_n'$ that is larger than $t_n$, and let $p > 0$. Then

$$h_n(t_n^*) < \frac{1}{p} + 2p \sqrt{\frac{2n - t_n}{t_n}} . \qquad \qquad \square$$

We now emphasize on the largest zero $t_n = T_n$ of $F_n$. By the results of § 2 the global maximum of $F_n$ is attained at the last zero $T_n^*$ of $F_n'$ which lies in the interval $(T_n, 2n)$, and is therefore the smallest zero of $F_n'$ after $T_n$. So Lemma 1 applies.

By a result of Bottema and Hahn (see [5], anffd [8], p. 228, last formula), the inequalityff

$$T_n > 2n - \frac{3}{2} - 8\sqrt{2}\sqrt{n-1} =: \tau_n \tag{43}$$

($n \geq 33$) holds for the last zero of $F_n$ (or $L_{n-1}^{(1)}$). As $\sqrt{\frac{2n-t}{t}}$ is strictly decreasing for $t \in (0, 2n)$, we have the inequality

$$\sqrt{\frac{2n - T_n}{T_n}} < \sqrt{\frac{2n - \tau_n}{\tau_n}} .$$

Puiseux expansion yields the asymptotic expression ($n \to \infty$)

$$\sqrt{\frac{2n - \tau_n}{\tau_n}} = 2\sqrt[4]{2}\,\frac{1}{n^{1/4}} + \frac{131}{16\sqrt[4]{2}}\,\frac{1}{n^{3/4}} + O\left(\frac{1}{n^{5/4}}\right) \ ,$$

especially is

$$\sqrt{\frac{2n - \tau_n}{\tau_n}} \sim \frac{1}{n^{1/4}} \ .$$

In our calculations the value $p$ was arbitrary, so we have the freedom to choose it properly. The asymptotics suggest the choice $p \sim n^{1/8}$. For any $a > 0$ we get therefore

$$h_n(T_n^*) < \frac{a}{n^{1/8}} + 2\,\frac{n^{1/8}}{a}\sqrt{\frac{2n - T_n}{T_n}} < \frac{a}{n^{1/8}} + 2\,\frac{n^{1/8}}{a}\sqrt{\frac{2n - \tau_n}{\tau_n}} \sim \frac{1}{n^{1/8}} < \frac{b}{n^{1/8}}$$

for some $b > 0$.

We choose the value $a = 2\sqrt[8]{2}$ (minimizing the leading term in the corresponding Puiseux expansion) and get the global estimate

$$h_n(T_n^*) < 2\sqrt[8]{2}\,\frac{1}{n^{1/8}} + \frac{\sqrt{\frac{3}{2} + 8\sqrt{2}\,\sqrt{n-1}}}{\sqrt[8]{2}\sqrt{2n - \frac{3}{2} - 8\sqrt{2}\,\sqrt{n-1}}}\,n^{1/8} \ .$$

Now we remember that $F_n$ takes its global maximum over $\mathbb{R}^+$ at the point $T_n^*$, and so does $h_n$. We therefore have for all $a > 0$, $n \in \mathbb{N}$ and $t > 0$ the inequality

$$|F_n(t)| < \sqrt{2\sqrt[8]{2}\,\frac{1}{n^{1/8}} + \frac{\sqrt{\frac{3}{2} + 8\sqrt{2}\,\sqrt{n-1}}}{\sqrt[8]{2}\sqrt{2n - \frac{3}{2} - 8\sqrt{2}\,\sqrt{n-1}}}\,n^{1/8}} = \frac{2\sqrt[16]{2}}{n^{1/16}} + \frac{131}{64\,2^{7/16}\,n^{9/16}} + O\left(\frac{1}{n^{17/16}}\right) . \quad (44)$$

We mention that we get a better asymptotic estimate if we use the sharper left hand inequality (20) instead of (43), set $\tau_n^* := -1 + 2n - C\,(2n - 1)$ ($C$ constant) leading to the asymptotic result

$$\sqrt{\frac{2n - \tau_n^*}{\tau_n^*}} \sim \frac{1}{n^{1/3}} \ ,$$

and therefore by the choice $p \sim n^{1/6}$ and the same procedure as above to the

**Theorem 4** For all $t \in \mathbb{R}^+$ we have the asymptotic inequality ($n \geq N$)

$$|F_n(t)| < \frac{c}{n^{1/12}}$$

for some $c > 0$, and in particular the limiting value

$$\lim_{n \to \infty} |F_n(t)| = 0 \ . \qquad \qquad \square$$

Obviously this theorem strengthens Theorem 3.

In principle (44) enables one to prove the statement

$$|F_n(t)| \leq \frac{2}{e}$$

for all $n \in \mathbb{N}$. Therefore one shows that the estimation function

$$E(n) = \sqrt{2\sqrt[8]{2}\frac{1}{n^{1/8}} + \frac{\sqrt{\frac{3}{2} + 8\sqrt{2}\sqrt{n-1}}}{\sqrt[8]{2}\sqrt{2n - \frac{3}{2} - 8\sqrt{2}\sqrt{n-1}}}n^{1/8}}$$

of (44) is decreasing, and as $E(17821075) > 2/e$ and $E(17821076) < 2/e$, it remains to prove the result for only a finite number of initial values.

The number of initial values, however, can be decisively fflowered using that by (40) $|F(t)| \leq 2/e$ whenever $T_n^*/n \geq \frac{e}{\sqrt{2}}$, especially if $T_n/n \geq \frac{e}{\sqrt{2}}$. From the Bottema-Hahn bound

$$\frac{T_n}{n} > \frac{\tau_n}{n} = 2 - \frac{3}{2n} - 8\sqrt{2}\frac{\sqrt{n-1}}{n} =: e(n)$$

we obtain first by the calculation

$$e'(n) = \frac{-16\sqrt{2} + 3\sqrt{n-1} + 8\sqrt{2}\,n}{2\,n^2\sqrt{n-1}}$$

that $e(n)$ is increasing for $n \geq 2$, and as $\lim\limits_{n \to \infty} e(n) = 2$ there is exactly one solution $n_0 \geq 2$ of the equation $e(n) = \frac{e}{\sqrt{2}}$, and $T_n/n > e(n) \geq e(n_0) = \frac{e}{\sqrt{2}}$ for $n > n_0$. A numerical calculation shows that $n_0 \approx 21138.7$ so that we are lead to the

**Theorem 5** The inequality $|F_n(t)| \leq 2/e$ is true for all $n \in \mathbb{N}$ and all $t > 0$ if it is true for $n \leq 21138$.

## 7 Estimates for the derivative

By (23) it follows that at the zeros of $F_n$ the derivative $F_n'$ satisfies the relation $|F_n'(t)| \leq 2$. This result holds for all $t \geq 0$ which can be seen as follows: Using (36) with $\alpha = -1$ we have

$$F_n'(t) = -e^{-t}\left(L_n(2t) + L_{n-1}(2t)\right), \tag{45}$$

and by an application of the Szegö result (28) it follows for $t > 0$

$$|F_n'(t)| < 2. \tag{46}$$

This shows that for all $n \in \mathbb{N}$ the derivatives $|F_n'|$ have their maximal value at the origin where $|F_n'(0)| = 2$ (see (11)). We note moreover, even though ffthe derivative $F := F_n'$ has a representation (45) similar to that of $F_n$ itself, it satisfies the much more complicated differential equation

$$\left(-4\,n^2 + 4\,n\,t - t^2\right)F(t) - 2\,n\,F'(t) + \left(-2\,n\,t + t^2\right)F''(t) = 0.$$

## 8   The functions $H_n$

In this section we collect the explicit inequalities that we deduced, and formulate a conjecture concerning the Bateman functions.

As the last point of inflection of the functions $F_n$ is at the point $t = 2n$ which increases with increasing $n$, it is reasonable to introduce the functions

$$H_n(t) := (-1)^n F_n(nt)$$

that have common absolute values with $F_n$ which, however, are attained at different points. The scale on the $t$-axis is here such that the last point of inflection lies at $t = 2$ for all functions $H_n$ ($n \in \mathbb{N}$), and $H_n$ is positive for $t \geq 2$. It is easy to deduce the differential equation

$$t\, H_n''(t) = n^2\, (t - 2)\, H_n(t) \tag{47}$$

satisfied by $H_n$.

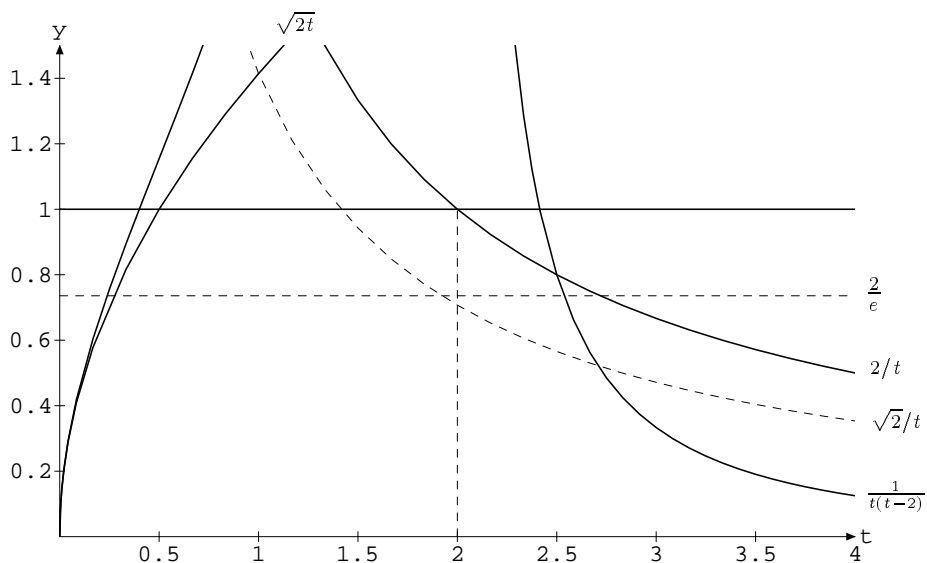

Figure 3: Estimates for the functions $H_n$

The inequalities that we deduced for $F_n$ read as follows for $H_n$: The trivial estimate (18) gives

$$|H_n(t)| \leq 1 \ ,$$

(16) yields ($n > 2$)

$$|H_n(t)| \leq \frac{2}{t} \ ,$$

the refinement (39) gives

$$|H_n(t_k)| \leq \frac{\sqrt{2}}{t_k}$$

for a critical point $t_k$ of $H_n$, (26) implies

$$|H_n(t)| < \sqrt{\frac{4t}{2 - t}} \ ,$$

and finally (34) yields

$$|H_n(t)| < \sqrt{2t} \; .$$

These estimates commonly do not depend on $n$. One more estimate will be added in the next section. Figure 3 shows them graphically.

Figures 4 and 5 show the graphs of the functions $H_n$ ($n = 1, \ldots, 20$). We conjecture that $H_n$ is strictly decreasing for increasing $n$ at the point $t = 2$. Note that by the result of Hahn (21) this is not true for any $t < 2$. In the next section we will show, however, that $\lim_{n \to \infty} H_n(t) = 0$ for each $t > 2$.

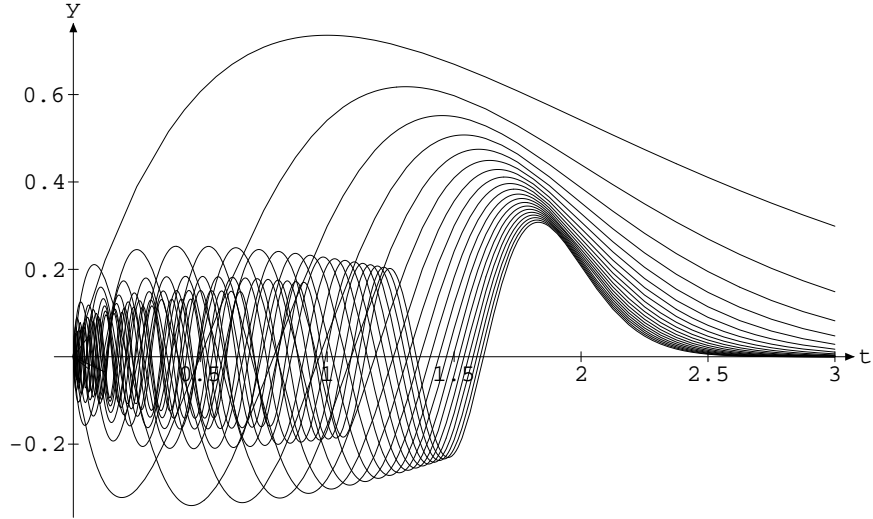

Figure 4: The functions $H_n$ ($n = 1, \ldots, 20$)

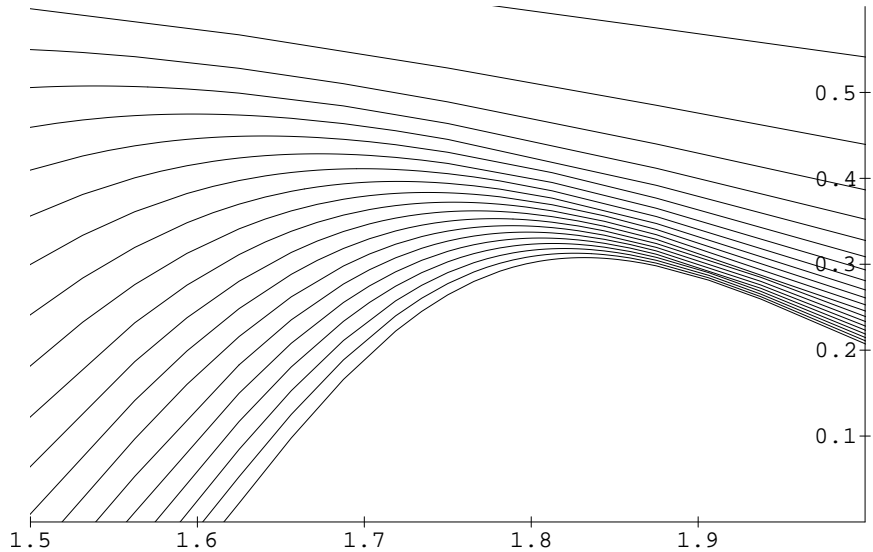

Figure 5: The functions $H_n$ ($n = 1, \ldots, 20$) in the interval $[1.5, 2]$

# 9 Estimates for large $t$

In this section we show that for all $t > 2$ the values $H_n(t)$ tend to 0 for $n \to \infty$.

The inequalities (46), and (17) correspond to the inequalities

$$|H_n'(t)| \leq 2n \ ,$$

and

$$|H_n'(t)| \leq \frac{n}{t} \qquad (48)$$

for the derivative of $H_n$, respectively.

By the differential equation (47) for $H_n$ we have

$$H_n(t) = \frac{t}{n^2 (t-2)} H_n''(t) \ . \qquad (49)$$

Let now $2 < t_1 < t_2$ be given. By definition $H_n$ is strictly positive in $[t_1, t_2]$, and by (49) $H_n''$ is strictly positive in $[t_1, t_2]$. Therefore

$$\int_{t_1}^{t_2} H_n(t)\, dt = \frac{1}{n^2} \int_{t_1}^{t_2} \frac{t}{t-2} H_n''(t)\, dt \ .$$

As the function $t/(t-2)$ is strictly decreasing, we have

$$\max_{t \in [t_1, t_2]} \frac{t}{t-2} = \frac{t_1}{t_1 - 2} \ ,$$

and therefore

$$\frac{1}{n^2} \int_{t_1}^{t_2} \frac{t}{t-2} H_n''(t)\, dt \leq \frac{1}{n^2} \frac{t_1}{t_1 - 2} \int_{t_1}^{t_2} H_n''(t)\, dt \ .$$

As $H_n''$ is positive in $[t_1, t_2]$, and as $\lim_{t \to \infty} H_n'(t) = 0$, it follows that $H_n'$ is negative and increasing, and therefore

$$\int_{t_1}^{t_2} H_n(t)\, dt \leq \frac{1}{n^2} \frac{t_1}{t_1 - 2} |H_n'(t_1) - H_n'(t_2)| \leq \frac{1}{n^2} \frac{t_1}{t_1 - 2} |H_n'(t_1)| \leq \frac{1}{n^2} \frac{t_1}{t_1 - 2} \frac{n}{t_1} \leq \frac{1}{n (t_1 - 2)}$$

where we used (48). For fixed $t_2 > 2$ we set now $t_1 := \frac{2 + t_2}{2}$, and with the integral mean value theorem we find $\tau \in [t_1, t_2]$ with

$$H(t_2) \leq H(\tau) = \frac{1}{t_2 - t_1} \int_{t_1}^{t_2} H_n(t)\, dt \leq \frac{1}{(t_2 - t_1)(t_1 - 2)} \frac{1}{n} = \frac{4}{(t_2 - 2)^2} \frac{1}{n} \ .$$

Another estimate for large $t$ which is independent of $n$, will be established now. Let again $2 < t_1 < t_2$. Then by (49)

$$
\begin{aligned}
|H_n^2(t_2) - H_n^2(t_1)| &= \left| \int_{t_1}^{t_2} (H_n^2(t))'\, dt \right| = \left| \int_{t_1}^{t_2} 2\, H_n'(t)\, H_n(t)\, dt \right| = \left| \int_{t_1}^{t_2} \frac{t}{n^2 (t-2)} H_n'(t)\, H_n''(t)\, dt \right| \\
&\leq \frac{t_1}{n^2 (t_1 - 2)} \int_{t_1}^{t_2} |2\, H_n'(t)\, H_n''(t)|\, dt = \frac{t_1}{n^2 (t_1 - 2)} |(H_n')^2(t_2) - (H_n')^2(t_1)| \ .
\end{aligned}
$$

We let now $t_2 \to \infty$, and get ($t := t_1$) using (48)

$$|H_n(t)|^2 \leq \frac{t}{n^2\,(t-2)}|(H_n')^2(t)| \leq \frac{1}{t\,(t-2)} \ .$$

This estimate improves the earlier ones for large $t$, see Figure 3.

Finally we show that uniformly with respecfft to $n$ the functions $H_n$ decrease faster than each negative power.

**Theorem 6** For each $k \in \mathbb{N}_0$ there is a constant $C_k > 0$ such that

$$|H_n(t)| \leq \frac{C_k}{t^k} \qquad (t > 0,\ n \in \mathbb{N})$$

that is independent of $n$.

*Proof:*  We prove the result by induction with respect to $k$. For $k = 0$ the statement is trivially true, see (18). Assume now the statement holds for some $k \in \mathbb{N}_0$, i. e.

$$|F_n(t)| \leq C_k\,\left(\frac{n}{t}\right)^k \ .$$

Then we get using (14)

$$
\begin{aligned}
|F_n(t)| &= \frac{1}{2t}\left|(n-1)\,(F_n(t)-F_{n-1}(t))+(n+1)\,(F_n(t)-F_{n+1}(t))\right| \\
&\leq \frac{2\,(n-1)\,C_k+(n+1)\,C_k+(n+1)\,C_k\,D_k}{2t}\,\left(\frac{n}{t}\right)^k \\
&\leq \frac{2\,(n-1)\,C_k\,D_k+(n+1)\,C_k\,D_k+(n+1)\,C_k\,D_k}{2t}\,\left(\frac{n}{t}\right)^k \\
&= 2\,C_k\,D_k\,\left(\frac{n}{t}\right)^{k+1} \ ,
\end{aligned}
$$

where we chose $D_k \geq 1$ such that $(n+1)^k \leq D_k\,b^k$ (the choice $D_k = 2^k$ does the job required as $(n+1)^k \leq (2n)^k \leq 2^k\,n^k$). This yields the result. $\qquad\square$